\documentclass[a4paper,10pt]{article}
\usepackage{mypackages}
\usepackage{mymacros}
\usepackage{mystyle}
\title{A totally bounded Uniformity on coarse metric Spaces}
\author{Elisa Hartmann}

\begin{document}

\maketitle

\begin{abstract}
This paper presents a new version of boundary on coarse spaces. The space of ends functor maps coarse metric spaces to uniform topological spaces and coarse maps to uniformly continuous maps.
\end{abstract}

\tableofcontents

\section{Introduction}

Coarse Geometry of metric spaces studies the large scale properties of a metric space. Meanwhile uniformity of metric spaces is about small scale properties.

Our purpose is to pursue a new version of duality between the coarse geometry of metric spaces and uniform spaces. We present a notion of boundary on coarse metric spaces which is a totally bounded separating uniform space. The methods are very basic and do not require any deep theory.

Note that the topology of metric spaces is well understood and there are a number of topological tools that can  be applied on coarse metric spaces which have not been used before. The new discovery may lead to new insight on the topic of coarse geometry.

\subsection{Background and related Theories}

There are quite number of notions for a boundary on a metric space. In this chapter we are going to discuss properties for three of them. The first paragraph is denoted to the \emph{Higson corona}, in the second paragraph the \emph{space of ends} is presented and in the last paragraph we study the \emph{Gromov boundary}.  

First we present the Higson corona. If $X$ is a proper metric space the \emph{Higson corona} $\nu X$ is the boundary of the Higson compactification $hX$ of $X$ which is a compact topological space that contains the underlying topological space of $X$ as a dense open subset. 

If $C(X)$ denotes the bounded continuous functions on $X$ then the so called Higson functions are a subset of $C(X)$. This subset determines a compactification which is called the Higson compactification. By a comment on~\cite[p.~31]{Roe2003} the Higson corona can be defined for any coarse space. The same does not work for the Higson compactification\footnote{for which the topology of $X$ needs to be locally compact which is given if the metric is proper.}. The~\cite[Proposition~2.41]{Roe2003} implies that the Higson corona is a covariant functor that sends coarse maps modulo closeness to continuous maps. Thus $\nu$ is a functor:
 \[
 \nu:\coarse\to \topology
 \]
 The topology of $\nu X$ has been studied in~\cite{Keesling1997}. It was shown in~\cite[Theorem~1]{Keesling1997} that for every $\sigma$-compact subset $A\s\nu X$ the closure $\bar A$ of $A$ in $\nu X$ is equivalent to the Stone-\v Cech compactification of $A$. The topology of $\nu X$ is quite complicated, especially if $X$ is a metric space. It has been noted in~\cite[Exercise~2.49]{Roe2003} that the topology of $\nu X$ for $X$ an unbounded proper metric space is never second countable. In~\cite[Theorem~1.1]{Dranishnikov1998} and~\cite[Theorem~7.2]{Dranishnikov2000} it was shown that if the asymptotic dimension $\asdim(X)$ of $X$ is finite then
\[
 \asdim(X)=\dim(\nu X)
\]
where the right side denotes the topological dimension of $\nu X$. Note that one direction of the proof uses the notion
of coarse covers\footnote{but under a different name.}.

Now we present the \emph{space of ends}. If $Y$ is a locally connected, connected and locally compact Hausdorff space then the space of ends $\Omega Y$ of $Y$ is the boundary of the Freudenthal compactification $\varepsilon Y$. It is totally disconnected and every other compactification of $Y$ that is totally disconnected factors uniquely through $\varepsilon Y$ by~\cite[Theorem~1]{Peschke1990}. The points of $\Omega Y$ are called \emph{endpoints} or \emph{ends}. 

Now~\cite[Theorem~5]{Peschke1990} shows that if $Y$ is a connected locally finite countable CW-complex every endpoint of $Y$ can be represented by a proper map 
\[
a:\R_+\to Y.
\]
Two proper maps $a_1,a_2:\R_+\to Y$ represent the same endpoint if they are connected by a proper homotopy. Denote by $\texttt{pTop}$ the category of topological spaces and proper continuous maps. Then the association $\Omega \cdot$ is a functor:
\[
\Omega:\texttt{pTop}\to \topology
\]
If $Y$ is a locally compact Hausdorff space then $\Omega Y$ can be constructed using a proximity relation which is a relation on the subsets of $Y$. See~\cite{Kron2009} for that one.

 This section studies the Gromov boundary. If $X$ is a proper Gromov hyperbolic metric space then the \emph{Gromov boundary} $\partial X$ consists of equivalence classes of \emph{sequences that converge to infinity} in $X$. The topology on $\partial X$ is generated by a basis of open neighborhoods. Loosely speaking two points on the boundary are close if the sequences that represent them stay close for a long time.
 
 By~\cite[Proposition~2.14]{Benakli2002} the topological spaces $\partial X$ and $\partial X\cup X$ are compact and by \cite[Theorem~2.1]{Benakli2002} the topology on $\partial X$ is metrizable. If $f:X\to Y$ is a \emph{quasi-isometry} between proper Gromov hyperbolic groups then it extends to a homeomorphism
\[
 \partial f:\partial X\to \partial Y
\]
by~\cite[Proposition~2.20]{Benakli2002}. In~\cite{Dydak2016} is studied a notion of morphisms for which the Gromov boundary is a functor: If $f:X\to Y$ is a \emph{visual function} between proper Gromov hyperbolic metric spaces then there is an induced map 
\[
 \partial f:\partial X\to \partial Y
\]
which is continuous by~\cite[Theorem~2.8]{Dydak2016}.

Now is there a notion of boundary on metric space which is both a functor on coarse spaces and coarse maps and has nice properties such as being Hausdorff and locally compact. As it turns out there is one such functor which is going to be designed in the course of this study. 

\subsection{Main Contributions}

 In the course of this study we will define a functor that associates to every coarse metric space a space at infinity which is a topological space. 

Based on the observation that twisted coarse cohomology with $\Z/2\Z$-coefficients of $\Z^n$ is the same as singular cohomology of $S^{n-1}$ with $\Z/2\Z$ coefficients we considered notions of boundary which reflect that observation.

In Definition~\ref{defn:closerelation} we define a relation on subsets of a metric space. As it turns out this relation is almost but not quite a proximity relation as noted in Remark~\ref{rem:proximityrelation}. The proof of Proposition~\ref{prop:starrefinement} uses that $X$ is a metric space, it does not work for general coarse spaces. That is why we restrict our study to metric spaces.

Note that while constructing the functor we presuppose which kind of spaces we want to distinguish. Indeed there is a certain class of metric spaces for which the local structure looks boring. The functor that we are going to define, the space of ends functor, is well suited for metric spaces that are coarsely proper coarsely geodesic. That class includes all Riemannian manifolds and finitely generated groups.

While the topology of the space at infinity is immediately defined using coarse covers there are two choices of points which are both solid: If $X$ is a metric space

\begin{itemize}
\item[(A):] the endpoints of $X$ are images of coarse maps $\Z_+\to X$ modulo finite Hausdorff distance or
\item[(B):] the points at infinity are subsets of $X$ modulo finite Hausdorff distance.
\end{itemize}

Note that choice B has been implemented in \cite{Grzegrzolka2018}. The space at infinity with choice A contains strictly less points than choice B. The~\cite[Proposition 93]{Hartmann2017a} guarantees that for choice A there exists at least one endpoint if the space $X$ is coarsely proper coarsely geodesic. The proof of~\cite[Proposition 93]{Hartmann2017a} is similar to the one of Königs Lemma in graph theory.

The space at infinity functor with choice B reflects isomorphisms by \cite[Proposition 2.18]{Kalantari2016} and the space at infinity functor with choice A is representable.

In the course of this article and in Definition~\ref{defn:endpoints} we use choice A, endpoints are images of coarse maps $\Z_+\to X$. Then we define the topology of the space of ends, $\soe(X)$, via surroundings of the diagonal in Definition~\ref{defn:spaceofends}. The uniformity on $\soe(X)$ is generated by a basis $(D_\ucover)_\ucover$ of entourages over coarse covers $\ucover$ of $X$. If $f:X\to Y$ is a coarse map then it induces a uniformly continuous map $\soe(f):\soe(X)\to \soe(X)$ between spaces of ends. That way the space of ends $\soe(\cdot)$ is a functor, we obtain the following result:

\begin{thma}
If $\mcoarse$ denotes the category of metric spaces and coarse maps modulo closeness and $\topology$ the category of 
topological spaces and continuous maps then $\soe$ is a functor
\[
 \soe:\mcoarse\to \topology.
\]
If $\uniformity$ denotes the category of uniform spaces and uniformly continuous maps then $\soe$ is a functor
\[
 \soe:\mcoarse\to \uniformity.
\]
\end{thma}

 It was nontrivial to show that a subspace in the domain category gives rise to a subspace in the image category. Proposition~\ref{prop:inclusionE} shows if $Z\s Y$ is a subspace then the inclusion $i:Z\to Y$ induces a uniform embedding $\soe(i):\soe(Z)\to \soe(Y)$.
 
 The functor $\soe(\cdot)$ preserves coproducts by Lemma~\ref{lem:preservescoprod}. The uniformity on $\soe(X)$ is totally bounded by Lemma~\ref{lem:totallybounded} and separating by Proposition~\ref{prop:hausdorff}.
\begin{thma}
If $X$ is a metric space then $\soe(X)$ is totally bounded and separating.
\end{thma}

We still lack a good study including the most basic properties of the new space of ends functor like compact and metrizable probably because the proofs are more difficult. 

\subsection{Outline}
 
 The outline of this paper is as follows:
 \begin{itemize}
  \item Chapter~\ref{sec:metric} introduces metric spaces and the coarse geometry of them. We recall basic notation used in this study.
  \item The Chapter~\ref{sec:tbu} creates the basic tool sets which are going to be used in this study.
  \item Finally Chapter~\ref{sec:definition} is about the space of ends functor. We give the definition and prove that the space of ends is a functor.
  \item Chapter~\ref{sec:props} studies a few basic properties. We prove that the space of ends is totally bounded and separated.
  \item The study closes with Chapter~\ref{sec:side} which sets the space of ends functor in context with other notions.
 \end{itemize}

\section{Metric Spaces}
\label{sec:metric}
\begin{defn}
 Let $(X,d)$ be a metric space. 
 \begin{itemize}
 \item Then the \emph{bounded coarse structure associated to $d$} on $X$ consists of those subsets $E\s X^2$ for which
 \[
  \sup_{(x,y)\in E}d(x,y)<\infty.
 \]
 We call an element of the coarse structure \emph{entourage}. 
\item The \emph{bounded cocoarse structure associated to $d$} on $X$ consists of those subsets $C\s X^2$ such that every sequence $(x_i,y_i)_i$ in $C$ is either bounded (which means both of the sequences $(x_i)_i$ and $(y_i)_i$ are bounded) or the set $d(x_i,y_i)_i$ is not bounded. 
 for $i\to \infty$. We call an element of the cocoarse structure coentourage.
 \item In what follows we assume the metric $d$ to be finite for every $(x,y)\in X^2$.
 \end{itemize}
\end{defn}

\begin{rem}
 Note that there is a more general notion of coarse spaces. By~\cite[Theorem~2.55]{Roe2003} a coarse structure on a coarse space $X$ is the bounded coarse structure associated to some metric $d$ on $X$ if and only if the coarse structure has a countable base. 
\end{rem}

\begin{rem}
 If $X$ is a metric space a subset $B\s X$ is \emph{bounded} if the set $B^2$ is an entourage in $X$.
\end{rem}

\begin{rem}
 Note the following duality:
 \begin{itemize}
  \item A subset $F\s X^2$ is an entourage if and only if for every coentourage $C\s X^2$ there is a bounded set $A\s X$ such that
\[
 F\cap C\s A^2.
\]
\item A subset $D\s X^2$ is a coentourage if and only if for every entourage $E\s X^2$ there is a bounded set $B$ such that
\[
 E\cap D\s B^2.
\]
 \end{itemize}
\end{rem}

\begin{defn}
 A map $f:X\to Y$ between metric spaces is called \emph{coarse} if
 \begin{itemize}
  \item $E\s X^2$ being an entourage implies that $\zz f E$ is an entourage \emph{(coarsely uniform)};
  \item and if $A\s Y$ is bounded then $\iip f A$ is bounded \emph{(coarsely proper)}.
 \end{itemize}
 Or equivalently
 \begin{itemize}
  \item $B\s X$ being bounded implies that $f(B)$ is bounded;
  \item and if $D\s Y^2$ is a coentourage then $\izp f D$ is a coentourage.
 \end{itemize}
 Two maps $f,g:X\to Y$ between metric spaces are called \emph{close} if
 \[
 f\times g(\Delta_X)
 \]
 is an entourage in $Y$. Here $\Delta_X$ denotes the diagonal in $X$.
\end{defn}

\begin{notat}
 A map $f:X\to Y$ between metric spaces is called
 \begin{itemize}
 \item \emph{coarsely surjective} if there is an entourage $E\s Y^2$ such that 
 \[
 E[\im f]=Y
 \]
 \item \emph{coarsely injective} if 
 \begin{enumerate}
  \item for every entourage $F\s Y^2$ the set $\izp f F$ is an entourage in $X$.
  \item or equivalently if for every coentourage $C\s X^2$ the set $\zzp f C$ is a coentourage in $Y$.
 \end{enumerate}
  \item two subsets $A,B\s X$ are called \emph{coarsely disjoint} if $A\times B$ is a coentourage.
\end{itemize}
\end{notat}

\begin{rem}
 We study metric spaces up to coarse equivalence. A coarse map $f:X\to Y$ is a \emph{coarse equivalence} if
 \begin{itemize}
  \item There is a coarse map $g:Y\to X$ such that $f\circ g$ is close to $id_Y$ and $g\circ f$ is close to $id_X$.
 \item or equivalently if $f$ is both coarsely injective and coarsely surjective.
 \end{itemize}
\end{rem}

\section{Totally Bounded Uniformity}
\label{sec:tbu}
\begin{defn}\name{close relation}
\label{defn:closerelation}
 Let $X$ be a metric space. Two subsets $A,B\s X$ are called \emph{close} if they are not coarsely disjoint. We write 
 \[
  A\close B.
 \]
Then $\close$ is a relation on the subsets of $X$.
\end{defn}

\begin{lem}
\label{lem:close}
 In every metric space $X$:
 \begin{enumerate}
  \item if $B$ is bounded, $B\notclose A$ for every $A\s X$
  \item $U\close V$ implies $V\close U$
  \item $U\close (V\cup W)$ if and only if $U\close V$ or $U\close W$
 \end{enumerate}
\end{lem}
\begin{proof}
\begin{enumerate}
 \item easy.
 \item easy.
 \item easy.
\end{enumerate}
\end{proof}

\begin{prop}
\label{prop:starrefinement}
Let $X$ be a metric space. Then for every subspaces $A,B\s X$ with $A\notclose B$ there are subsets $C,D\s X$ such that $C\cap D=\emptyset$ and $A\notclose (X\ohne C)$, $B\notclose X\ohne D$.
\end{prop}
\begin{proof}
Note this is the same as~\cite[Proposition~4.5]{Kalantari2016} where the same statement was proven in a similar fashion. Let $E_1\s E_2\s \cdots$ be a symmetric basis for the coarse structure of $X$. Then for every $x\in A^c\cap B^c$ there is a least number $n_1(x)$ such that $x\in E_{n_1(x)}[A]$ and a least number $n_2(x)$ such that $x\in E_{n_2(x)}[B]$. Define:
\[
 V_1=\{x\in A^c\cap B^c:n_1(x)\le n_2(x)\}
\]
and
\[
 V_2=A^c\cap B^c\ohne V_1.
\]
Now for every $n$ 
\[
 E_n[V_1]\cap B\s E_{2n}[A]\cap B 
\]
because for every $x\in V_1$, if $x\in E_n[B]$ then $x\in E_n[A]$. Now define
\[
 C=A\cup V_1
\]
and
\[
 D=B\cup V_2.
\]
\end{proof}
\begin{rem}
\label{rem:proximityrelation}
Compare $\close$ with the notion of proximity relation~\cite[chapter~40, pp.~266]{Willard1970}. By Lemma \ref{lem:close} and Proposition \ref{prop:starrefinement} the close relation satisfies \cite[P-1),P-3)-P-5) of Definition 40.1]{Willard1970} but not P-2).
\end{rem}

\begin{rem}
 \label{rem:coarsecontinuous}
 If $f:X\to Y$ is a coarse map then whenever $A\close B$ in $X$ then $f(A)\close f(B)$ in $Y$.
\end{rem}

We recall~\cite[Definition~45]{Hartmann2017a}:
\begin{defn}\name{coarse cover}
 If $X$ is a metric space and $U\s X$ a subset a finite family of subsets $U_1,\ldots,U_n\s U$ is said to \emph{coarsely cover} $U$ if
\[
 U^2\cap (\bigcup_i U_i^2)^c
\]
is a coentourage in $X$.
\end{defn}

\begin{rem}
 Note that coarse covers determine a Grothendieck topology on $X$. If $f:X\to Y$ is a coarse map between metric spaces and $(V_i)_i$ a coarse cover of $V\s Y$ then $(\iip f{V_i})_i$ is a coarse cover of $\iip f V\s X$.
\end{rem}

\begin{lem}
 Let $X$ be a metric space. A finite family $\ucover=\{U_\alpha:\alpha\in A\}$ is a coarse cover if and only if there is a finite cover $\vcover=\{V_\alpha:\alpha\in A\}$ of $X$ as a set such that $V_\alpha\notclose U_\alpha^c$ for every $\alpha$. 
\end{lem}
\begin{proof}
Suppose $\ucover$ is a coarse cover of $X$. We proceed by induction on the index of $\ucover$:
 \begin{itemize}
  \item $n=1$: a subset $U$ coarsely covers $X$ if and only if $U^c$ is bounded if and only if $U^c\notclose X$.
  \item two subsets $U_1,U_2$ coarsely cover $X$ if and only if $U_1^c\notclose U_2^c$. Now by Proposition~\ref{prop:starrefinement} there are $C,D\s X$ with $C^c\cup D^c=X$ and $C^c\notclose U_1^c$ and $D^c\notclose U_2^c$. Define $V_1=C^c,V_2=D^c$.
\item $n+1\to n+2$: Subsets $U,V,U_1,\ldots,U_n$ coarsely cover $X$ if and only if $U,V$ coarsely cover $U\cup V$ and $U\cup V,U_1,\ldots, U_n$ coarsely cover $X$ at the same time.

Suppose $U,V,U_1,\ldots,U_n$ coarsely cover $X$. By induction hypothesis there is a cover of sets $V_1',V_2'$ of 
$U\cup V$ such that $V_1'\notclose U^c\cap V$ and $V_2'\notclose V^c\cap U$ and there is a cover of sets $W,V_1,\ldots, V_n$ such 
that $W\notclose (U\cup V)^c$ and $V_i\notclose U_i^c$ for every $i$. Then $V_1'\cap W\notclose U^c,V_2'\cap W\notclose V^c$. Now
\[
B:=(U\cup V)^c\cap W
\]
is bounded.
Then 
\[
V_1'\cap W,V_2'\cap W,V_1,\ldots,V_n\cup B
\]
is a finite cover of $X$ with the desired properties.
\end{itemize}
Suppose $(V_\alpha)_\alpha$ cover $X$ as sets and $V_\alpha\notclose U_\alpha^c$ for every $\alpha$. Let $E\s X^2$ be an entourage. Then $E[U_\alpha^c]\cap V_\alpha$ is bounded for every $\alpha$. Then
\f{
\bigcap_\alpha E[U_\alpha^c]
&=\bigcap_\alpha E[U_\alpha^c]\cap(\bigcup_\alpha V_\alpha)\\
&=\bigcup_\alpha (V_\alpha\cap\bigcap_\beta E[U_\beta^c])
}
is bounded. By \cite[Proposition 54]{Hartmann2017a} we can conclude that $(U_\alpha)_\alpha$ is a coarse cover of $X$.
\end{proof}

\begin{rem}
The~\cite[Theorem~40.15]{Willard1970} states that every proximity relation on a set is induced by some totally bounded 
uniformity on it. Note that a coarse cover on a metric space $X$ does not precisely need to cover $X$ as a set. Except 
for that the collection of all coarse covers of a metric space satisfies \cite[a),b) of Theorem 36.2]{Willard1970}. We 
can compare coarse covers of $X$ with a base for a totally bounded uniformity on $X$: the collection of all sets 
$\bigcup_i U_i^2$ for $(U_i)_i$ a coarse cover satisfies \cite[b)-e) of Definition 35.2]{Willard1970} but not a). Note 
that by \cite[Definition 39.7]{Willard1970} a diagonal uniformity is totally bounded if it has a base consisting of 
finite covers.
\end{rem}

\begin{lem}\name{separation cover}
\label{lem:separationcover}
 If $U_1,U_2$ coarsely cover a metric space $X$ (or equivalently if $U_1^c,U_2^c$ are coarsely disjoint) then there exists a coarse cover $V_1,V_2$ of $X$ such that $V_1\notclose U_1^c$ and $V_2\notclose U_2^c$.
\end{lem}
\begin{proof}
By Proposition \ref{prop:starrefinement} there are subsets $C,D\s X$ such that $C\cap D=\emptyset$, $U_1^c\notclose C^c$ and $U_2^c\notclose D^c$. Thus $U_1,C$ is a coarse cover of $X$ such that $C\notclose U_2^c$. 

By Proposition \ref{prop:starrefinement} there are subsets $A,B\s X$ such that $A\cap B=\emptyset$, $A^c\notclose U_1^c$ and $B^c\notclose C^c$. Then $B,C$ are a coarse cover of $X$ such that $B\notclose U_1^c$.

Then $V_1=B$ and $V_2=C$ have the desired properties.
 \end{proof}

\begin{notat}\name{coarse star refinement}
 Let $\ucover=(U_i)_{i\in I}$ be a coarse cover of a metric space $X$.
 \begin{enumerate}
  \item If $S\s X$ is a subset then
  \[
   \coarsestar S \ucover=\bigcup \{U_i:S\close U_i\}
  \]
 is called the \emph{coarse star of $S$}.
 \item Let $\vcover=(V_j)_{j\in J}$ be a coarse cover of $X$. Then $\vcover$ is a \emph{coarse barycentric refinement of $\ucover$} if for every $j_1,\ldots,j_k\in J$ and entourage $E\s X^2$ with
 \[
  \bigcap_k E[V_{i_k}]
 \]
 not bounded there is some $i\in I$ and an entourage $F\s X^2$ such that
\[
\bigcup_k V_{i_k}\s F[U_i].
 \]
 \item A coarse cover $\vcover=(V_j)_{j\in J}$ of $X$ is called a \emph{coarse star refinement of $\ucover$} if for 
every $j\in J$ there is some $i\in I$ and entourage $E\s X^2$ such that
\[
 \coarsestar{V_j}\vcover \s E[U_i].
\]
 \end{enumerate}
\end{notat}

\begin{lem}
 \label{lem:coarsestarprops}
 If $\vcover=(V_i)_i$ is a coarse star refinement of a coarse cover $\ucover=(U_i)_i$ of a metric space $X$ then
\begin{itemize}
 \item if $S\s X$ is a subset then there is an entourage $E\s X^2$ such that
 \[
  \coarsestar{\coarsestar S \vcover}\vcover\s E[\coarsestar S \ucover];
 \]
 \item if $f:X\to Y$ is a coarse map between metric spaces, $(U_i)_i$ a coarse cover of $Y$ and $S\s X$ a subset then  \[
  f(\coarsestar S{\iip f \ucover})\s \coarsestar {f(S)} \ucover.
 \]
\end{itemize}
\end{lem}
\begin{proof}
 \begin{itemize}
\item Suppose $E\s X^2$ is an entourage such that for every $V_j$ there is an $U_i$ such that $\coarsestar {V_j} \vcover \s E[U_i]$. Note that $S\close V_j$ implies $S\close U_i$ in that case. Then
\f{
 \coarsestar {\coarsestar S \vcover} \vcover
 &= \coarsestar{\bigcup \{V_i: V_i\close S\}} \vcover\\
 &=\bigcup_{V_i\close S} \coarsestar {V_i} \vcover\\
 &\s \bigcup_{S\close U_j} E[U_j]\\
 &=E[\coarsestar S \ucover].
}
 \item 
 \f{
 f(\coarsestar S {\iip f \ucover})
&=\bigcup \{f\circ \iip f{U_i}: S\close \iip f {U_i}\}\\
&\s \bigcup \{f\circ \iip f {U_i}: f(S)\close f\circ \iip f{U_i}\}\\
&\s \bigcup \{U_i:f(S)\close U_i\}\\
&=\coarsestar {f(S)} \ucover.
 }
 \end{itemize}
\end{proof}

\begin{lem}
\label{lem:starrefinement}
 If $\ucover$ is a coarse cover of a metric space $X$ then there exists a coarse cover $\vcover$ of $X$ that coarsely star refines $\ucover$.
\end{lem}
\begin{proof}
There are three steps:
 
   If $\vcover=(V_j)_j$ is a coarse barycentric refinement of $\ucover$ and $\wcover=(W_k)_k$ is a coarse barycentric refinement of $\vcover$ then $\wcover$ is a coarse star refinement of $\ucover$: 
  
   fix $W_k$ and denote $J=\{j:W_k\close W_j\}$.
  
   Then for every $j\in J$ there is some $V_j$ and entourage $E_j\s X^2$ such that $W_k\cup W_j\s E_j[V_j]$.
  
   Define $E=\bigcup_j E_j$. Then $\bigcap_j E[V_j]\z W_k$. Thus there is some $U_i$ and entourage $F\s X^2$ such that $\bigcup_jV_j\s F[U_i]$. 
  
  For every $j\in J$:
   \f{
   W_j&\s E[V_j]\\
    &\s E\circ F[U_i].
   }
   Thus $\coarsestar {W_k} \wcover \s E\circ F[U_i]$.
  
   We show there is a coarse barycentric refinement $\vcover=(V_i)_i$ of $\ucover$: First we show if $U_1,U_2$ is a coarse cover of $X$ then there is a coarse barycentric refinement $V_1,V_2,V_3$ of $U_1,U_2$: 

   By Lemma \ref{lem:separationcover} there is a coarse cover $W_1,W_2$ of $X$ such that $W_1\notclose U_1^c$ and $W_2\notclose U_2^c$.
 
 Then $ W_1^c, U_1$ and $W_2^c, U_2$ are coarse covers of $X$.
 
 By Proposition \ref{prop:starrefinement} there are $C,D\s X$ such that $C\cap D=\emptyset$, $D^c\notclose U_2^c$, $C^c\notclose W_2$.
 
 Also there are $A,B\s X$ such that $A\cap B=\emptyset$, $A^c\notclose U_1^c$, $B^c\notclose W_1$.
 
 Then
 \[
  V_1=W_1,V_2=C\cap B,V_3=W_2
 \]
has the desired properties:

 $(V_i)_i$ is a coarse cover:  

 Note that by $B^c\notclose W_1$ and $W_1^c\notclose W_2^c$ the sets $W_2,B$ are a coarse cover of $X$.

 Note that by $C^c\notclose W_2$ and $W_1^c\notclose W_2^c$ the sets $W_1,C$ coarsely cover $X$.

 Note that $(W_1\cap W_2)\notclose (C^c\cup B^c)$. Then, combining items i,ii, we get that
\[
 W_1,W_2, B\cap C
\]
is a coarse cover as required.

 $(V_i)_i$ is a coarse barycentric refinement of $U_1,U_2$: 

 There is some entourage $E\s X^2$ such that $V_1\cup V_2\s E[U_1]$: For $W_1$ we use that $W_1\notclose U_1^c$. For $C\cap B$ we use that $A^c\notclose U_1^c$ and $B\s A^c$.  

 There is an entourage $E\s X^2$ such that $V_2\cup V_3 \s E[U_2]$: For $W_2$ we use that $W_2\notclose U_2^c$. For $C\cap B$ we use that $D^c\notclose U_2^c$ and $C\s D^c$.

 $V_1\notclose V_3$: We use $W_1^c\notclose W_2^c$.

  Now we show the general case: Suppose $U_i\s X$ are subsets such that $\ucover=(U_i)_i$. We show there is a coarse barycentric refinement $\vcover$ of $\ucover$.

 For every $i$ the sets $U_i,\bigcup_{j\not=i}U_j$ coarsely cover $X$. By Lemma \ref{lem:separationcover} there 
are subsets $W_1^i,W_2^i$ that coarsely cover $X$ such that $W_1^i\notclose U_i^c$ and $W_2^i\notclose 
(\bigcup_{j\not=i} U_j)^c$.

Then there is a coarse barycentric refinement $V^i_1,V^i_2,V^i_3$ of $W_1^i,W_2^i$ for every $i$.

Then we define 
\[
 \vcover:=(\bigcap_iV^i_{\sigma(i)})_\sigma,
\]
here $\sigma(i)\in\{1,2,3\}$ is all possible permutations.

We show $\vcover$ is a coarse cover of $X$ that is a coarse barycentric refinement of $\ucover$:

$\vcover$ is a coarse cover: by design.

$\vcover$ is a coarse barycentric refinement of $\ucover$: Suppose there is an entourage $E\s X^2$ and a subindex $(\sigma_k)_k$ such that $\bigcap_{\sigma_k} E[\bigcap_{\sigma_k}V_{\sigma_k(i)}^i]$ is not bounded. Then
\[
 \bigcap_{i,\sigma_k}E[V_{\sigma_k(i)}^i]
\]
is not bounded. Then there is an entourage $F\s X^2$ such that for every $i$:
\[
 \bigcup_{\sigma_k} V^i_{\sigma_k(i)}\s F[W^i_{l_i}]
\]
where $l_i$ is one of $1,2$. Then 
\f{
 \bigcup_{\sigma_k}\bigcap_i V^i_{\sigma_k(i)}
 &\s \bigcap_i \bigcup_{\sigma_k} V^i_{\sigma_k(i)}\\
 &\s \bigcap_i F[W^i_{l_i}]
}
if $l_i=1$ for one $i$ then we are done. Otherwise
\[
 \bigcup_{\sigma_k}\bigcap_i V^i_{\sigma_k(i)}\s \bigcap_i F[W^i_2]
\]
and $F[W_2^i]\notclose (\bigcup_{j\not=i} U_j)^c$ implies
\[
 \bigcap_i F[W_2^i] \notclose \bigcup_i (\bigcup_{j\not=i} U_j)^c
\]
which implies that $\bigcap_i F[W_2^i]$ is bounded, a contradiction.
\end{proof}

\section{Definition}
\label{sec:definition}
We introduce the space of ends of a coarse space which is a functor $\soe$ from the category of coarse metric spaces to the category of uniform spaces. 
\label{subsec:soe:def}
\begin{defn}\name{endpoint}
\label{defn:endpoints}
Let $X$ be a metric space,
\begin{itemize}
 \item two coarse maps $\phi,\psi:\Z_+\to X$ are said to represent the same \emph{endpoint in $X$} if there is an entourage $E\s X^2$ such that 
\[
  E[\psi(\Z_+)]=\phi(\Z_+).
\]
\item if $\ucover=(U_i)_i$ is a coarse cover of $X$ and $p,q$ are two endpoints in $X$ which are represented by $\phi,\psi:\Z_+\to X$. Then $q$ is said to be in a \emph{$\ucover-$neighborhood of $p$}, denoted $q\in\ucover[p]$, if there is an entourage $E\s X^2$ such that
 \[
 E[\coarsestar {\phi(\Z_+)}\ucover]\z \psi(\Z_+)
\]
and
 \[
 E[\coarsestar {\psi(\Z_+)}\ucover]\z \phi(\Z_+).
\]
\end{itemize}
\end{defn}

\begin{lem}
 If $\vcover\le \ucover$ is a refinement of a coarse cover of a metric space $X$ then for every two endpoints $p,q$ of $X$ the relation $q\in \vcover[p]$ implies the relation $q\in \ucover[p]$.
\end{lem}
\begin{proof}
 Suppose $\vcover=(V_i)_i$ and $\ucover=(U_i)_i$. If $p,q$ are represented by $\phi,\psi:\Z_+\to X$ then
 \f{
 \coarsestar{\phi(\Z_+)}\vcover &=\bigcup\{V_i:\phi(\Z_+)\close V_i\}\\
 &\s \bigcup\{U_i:\phi(\Z_+)\close U_i\}\\
 &=\coarsestar{\phi(\Z_+)}\ucover
 }
 in the same way $\coarsestar{\psi(\Z_+)}\vcover\s \coarsestar {\psi(\Z_+)}\ucover$. Then if $q\in\vcover[p]$ there is some entourage $E\s X^2$ such that
 \f{
 \psi(\Z_+)&\s E[\coarsestar{\phi(\Z_+)}\vcover]\\
 &\s E[\coarsestar{\phi(\Z_+)}\ucover]
 }
 and $\phi(\Z_+)\s E[\coarsestar{\psi(\Z_+)}\ucover]$. Thus $q\in\ucover[p]$. 
\end{proof}

\begin{defn}\name{space of ends}
\label{defn:spaceofends}
Let $X$ be a metric space. As a set the \emph{space of ends} $\soe(X)$ of $X$ consists of the endpoints in $X$. A subset $U\s \soe(X)$ is open if for every $p\in U$ there is a coarse cover $\ucover$ of $X$ such that 
\[
\ucover[p]\s U.
\]
This defines a topology on $\soe(X)$.
\end{defn}

\begin{rem}
 The topology on the set of endpoints $\soe(X)$ is generated by a uniformity: If $\ucover$ is a coarse cover of $X$ then 
 \[
  D_\ucover=\{(p,q):q\in \ucover[p]\}
 \]
is the \emph{entourage associated to $\ucover$}. Then $(D_\ucover)_\ucover$ over coarse covers $\ucover$ of $X$ are a base for a diagonal uniformity on $\soe(X)$.
\end{rem}

\begin{lem}
 If $X$ is a metric space then $\soe(X)$ is indeed a uniform space. Coarse covers of $X$ give rise to a base for the uniform structure.
\end{lem}
\begin{proof}
 We check that $(D_\ucover)_\ucover$ over coarse covers are a base for a uniformity on $\soe(X)$:
 \begin{enumerate}
  \item If $\ucover$ is a coarse cover of $X$ then $\Delta\s D_\ucover$, where $\Delta=\{(p,p):p\in E(X)\}$: $p\in \ucover[p]$.
  \item If $\ucover,\vcover$ are coarse covers of $X$ then $D_\ucover\cap D_\vcover$ is an entourage: Suppose $\ucover=(U_i)_i,\vcover=(V_i)_i$ then define
  \[
   \ucover\cap \vcover:=(U_i\cap V_j)_{ij}.
  \]
 Suppose $p,q$ are represented by $\phi,\psi:\Z_+\to X$. Then $q\in(\ucover\cap\vcover)[p]$ implies there is an entourage $E\s X^2$ such that
 \f{
 \psi(\Z_+)&\s E[\coarsestar{\phi(\Z_+)}{\ucover\cap \vcover}]\\
  &=E[\bigcup\{U_i\cap V_i:\phi(\Z_+)\close U_i\cap V_j\}\\
  &\s \bigcup_{U_i\cap V_j\close\phi(\Z_+)}E[U_i]\cap E[V_j]\\
  &\s (\bigcup_{U_i\close \phi(\Z_+)}E[U_i])\cap (\bigcup_{V_i\close \phi(\Z_+)}E[V_i])\\
  &=E[\coarsestar{\phi(\Z_+)}\ucover]\cap E[\coarsestar{\phi(\Z_+)}\vcover].
 }
Thus $q\in \ucover[p]\cap\vcover[p]$. This way we have proven:
 \[
  D_{\ucover\cap\vcover}\s D_\ucover\cap D_\vcover.
 \]
  \item If $\ucover$ is a coarse cover of $X$ then there is a coarse cover $\vcover$ of $X$ such that $D_\vcover\circ D_\vcover\s D_\ucover$: By Lemma~\ref{lem:starrefinement} there is a coarse star refinement $\vcover$ of $\ucover$. And by Lemma~\ref{lem:convergence1} item 2 the uniform cover $(\vcover[p])_p$ star refines the uniform cover $(\ucover[p])_p$ thus the result.
  \item If $\ucover$ is a coarse cover then $D_\ucover=D_\ucover^{-1}$.
 \end{enumerate}
A subset $D\s E(X)^2$ is an entourage of the uniform structure of $\soe(X)$ if there is a coarse cover $\ucover$ of $X$ such that 
\[
 D_\ucover\s D.
\]
\end{proof}

\begin{thm}
\label{thm:mapsE}
 If $f:X\to Y$ is a coarse map between metric spaces then the induced map
 \f{
  \soe(f):\soe(X)&\to \soe(Y)\\
  [\varphi]&\mapsto [f\circ \varphi]
 }
 is a continuous map between topological spaces.
\end{thm}
\begin{proof}
 We show $\soe(f)$ is well defined: if $\phi,\psi:\Z_+\to X$ represent the same endpoint in $X$ then there is some entourage $E\s X^2$ such that $E[\psi(\Z_+)]= \phi(\Z_+)$. But then
\f{
f^2(E)[f\circ\psi(\Z_+)]&\z f(E[\psi(\Z_+)])\\
& = f\circ\phi(\Z_+).
}
Thus $f\circ \phi,f\circ \psi$ represent the same endpoint in $Y$.

We show $\soe(f)$ continuous: 
For that we show that the inverse image of an open set is an open set. 

Let $U\s \soe(Y)$ be open and $p\in \soe(f)^{-1}(U)$ be a point. Suppose that $p$ is represented by a coarse map $\phi:\Z_+\to X$. Then $f\circ\phi$ represents $\soe(f)(p)\in U$. Now there is a coarse cover $\ucover=(U_i)_i$ of $Y$ such that $\ucover[\soe(f)(p)]\s U$. Then $\iip f \ucover=(\iip f {U_i})_i$ is a coarse cover of $X$.

If $q\in\iip f \ucover[p]$ we show that $\soe(f)(q)\in\ucover[\soe(f)(p)]$: Suppose that $q$ is represented by a coarse map $\psi:\Z_+\to X$. Then there is some entourage $F\s X^2$ such that
\[
 F[\coarsestar{\phi(\Z_+)}{\iip f \ucover}]\z \psi(\Z_+)
\]
and
\[
 F[\coarsestar{\psi(\Z_+)} {\iip f \ucover}]\z \phi(\Z_+).
\]
By Lemma \ref{lem:coarsestarprops}:
\f{
 f^2(F)[\coarsestar{f\circ\phi(\Z_+)}\ucover]
 &\z f^2(F)[f(\coarsestar{\phi(\Z_+)}{\iip f \ucover})])\\
  &\z f(F[\coarsestar{\phi(\Z_+)}{\iip f \ucover}])\\
  &\z f\circ\psi(\Z_+)
}
and
\f{
 f^2(F)[\coarsestar{f\circ\psi(\Z_+)}\ucover]
 &\z f^2(F)[f(\coarsestar{\psi(\Z_+)}{\iip f \ucover})])\\
  &\z f(F[\coarsestar{\psi(\Z_+)}{\iip f \ucover}])\\
  &\z f\circ\phi(\Z_+).
}
Now $f\circ \psi$ represents $\soe(f)(q)$ which by the above is in $\ucover[\soe(f)(p)]$.
\end{proof}

\begin{rem}
 The proof of Theorem \ref{thm:mapsE} uses the following: if $f:X\to Y$ is a coarse map and $D_\ucover$ the entourage of $\soe(Y)$ associated to a coarse cover $\ucover$ of $Y$ then there is an entourage $D_{\iip f \ucover}$ of $\soe(X)$ associated to the coarse cover $\iip f \ucover$ of $X$ such that $(p,q)\in D_{\iip f \ucover}$ implies $(\soe(f)(p),\soe(f)(q))\in D_\ucover$. Thus $\soe(f)$ is a uniformly continuous map between uniform spaces $\soe(X)$ and $\soe(Y)$.
\end{rem}

\begin{lem}
 If two coarse maps $f,g:X\to Y$ are close then $\soe(f)=\soe(g)$.
\end{lem}
\begin{proof}
 Let $p\in \soe(X)$ be a point that is represented by $\varphi$. Now $f,g$ are close thus $H:=f\times g(\Delta_X)$ is an entourage. But then
\[
 H[g\circ\varphi(\Z_+)]\z f\circ\varphi(\Z_+)
\]
thus $\soe(f)(p)=\soe(g)(p)$.
\end{proof}

\begin{cor}
\label{cor:isoE}
 If $f$ is a coarse equivalence then $\soe(f)$ is a homeomorphism between topological spaces $\soe(X)$ and $\soe(Y)$. In fact $\soe(f)$ is a uniform isomorphism between uniform spaces $\soe(X)$ and $\soe(Y)$.
\end{cor}

\begin{cor}
 If $\mcoarse$ denotes the category of metric spaces and coarse maps modulo closeness and $\topology$ the category of 
topological spaces and continuous maps then $\soe$ is a functor
\[
 \soe:\mcoarse\to \topology.
\]
If $\uniformity$ denotes the category of uniform spaces and uniformly continuous maps then $\soe$ is a functor
\[
 \soe:\mcoarse\to \uniformity.
\]
\end{cor}

\begin{ex}
 $\soe(\Z_+)$ is a point.
\end{ex}

\section{Properties}
\label{sec:props}
\begin{lem}
 \label{lem:convergence1}
 If $X$ is a metric space
 \begin{itemize}
 \item and $\vcover$ is a coarse star refinement of a coarse cover $\ucover$ of $X$ then $q\in \vcover[p]$ and 
$r\in\vcover[q]$ implies $r\in\ucover[p]$.
 \item if $\vcover$ coarsely star refines $\ucover$ then $(\vcover[p])_p$ star refines $(\ucover[p])_p$
 \end{itemize}
\end{lem}
\begin{proof}
 \begin{itemize}
  \item Suppose $p$ is represented by $\phi:\Z_+\to X$, $q$ is represented by $\psi:\Z_+\to X$ and $r$ is represented by $\rho:\Z_+\to X$. Then $E[\coarsestar {\phi(\Z_+)} \vcover\z \psi(\Z_+)$ and $E[\coarsestar {\psi(\Z_+)} \vcover\z \phi(\Z_+)$, $E[\coarsestar {\rho(\Z_+)} \vcover\z \psi(\Z_+)$ and $E[\coarsestar {\psi(\Z_+)} \vcover\z \rho(\Z_+)$. By Lemma~\ref{lem:coarsestarprops} there is an entourage $F\s X^2$ such that $\coarsestar {\coarsestar {\phi(\Z_+)}\vcover} \vcover \s F[\coarsestar {\phi(\Z_+)}\ucover]$. Then 
  \f{
   E^{\circ 2}\circ F[\coarsestar {\phi(\Z_+)}\ucover]
  &\z E^{\circ 2}[\coarsestar {\coarsestar {\phi(\Z_+)}\vcover} \vcover]\\
   &\z E[\coarsestar {\psi(\Z_+)}\vcover]\\
  &\z \rho(\Z_+)
  }
 the other direction works the same way.
  \item Fix $p\in \soe(X)$. Then
\[
 st(\vcover[p],(\vcover[p])_p)\s \ucover[p]
\]
because if $q\in \vcover[p]$ and $q\in \vcover[r]$ then $r\in \ucover[p]$ by Item 1
 \end{itemize}
\end{proof}

\begin{prop}
\label{prop:inclusionE}
 If $i:Z\to Y$ is an inclusion of metric spaces then $\soe(i):\soe(Z)\to \soe(Y)$ is a uniform embedding.
\end{prop}
\begin{proof}
That $\soe(i)$ is injective is easy to see.

Define a map 
\f{
\Phi:\soe(i)(\soe(Z))&\to \soe(Z)\\
\soe(i)(p)&\mapsto p.
}
We show $\Phi$ is a uniformly continuous map:

If $\ucover=(U_i)_i$ is a coarse cover of $Z$ we show there is a coarse cover $\vcover$ of $Y$ such that for every $p,q\in \soe(Z)$: the relation $\soe(i)(q)\in \vcover[\soe(i)(p)]$ implies $q\in \ucover[p]$.

Note that for every $i$ the sets $U_i^c\notclose (\bigcup_{j\not=i}U_j)^c$ are coarsely disjoint in $Y$. By Lemma 
\ref{lem:separationcover} there are subsets $W_1^i,W_2^i\s Y$ that coarsely cover $Y$ and $W^i_1\notclose U_i^c$ and $W^i_2\notclose (\bigcup_{j\not=i}U_j)^c$. Now define
\[
 \vcover:=(\bigcap_i W^i_{\sigma(i)})_\sigma
\]
with $\sigma(i)\in \{1,2\}$ all possible permutations. Note that there is some entourage $E\s Z^2$ such that for every 
$\sigma$ there is some $U_j$ such that
\[
 \bigcap_i W^i_{\sigma(i)}\cap Z\s E[U_j].
\]

Let $p,q\in \soe(Z)$ such that $\soe(i)(q)\in \vcover[\soe(i)(p)]$. Suppose $p,q$ are represented by $\phi,\psi:\Z_+\to 
Z$. Then there is some entourage $F\s Y^2$ such that 
\[
 F[\coarsestar {\phi(\Z_+)}\vcover]\z \psi(\Z_+)
\]
and
\[
 F[\coarsestar {\psi(\Z_+)}\vcover]\z \phi(\Z_+).
\]
Then
\f{
F\circ E[\coarsestar{\phi(\Z_+)}\ucover]
&\z F[\coarsestar{\phi(\Z_+)}\vcover\cap Z]\\
&\z \psi(\Z_+).
}
The other direction works the same way.

Then $q\in \ucover[p]$ as we wanted to show.
\end{proof}

\begin{rem}
By Proposition~\ref{prop:inclusionE} and Corollary~\ref{cor:isoE} every coarsely injective coarse map $f:X\to Y$ induces a uniform embedding. We identify $\soe(X)$ with its image $\soe(f)(\soe(X))$ in $\soe(Y)$.
\end{rem}

\begin{ex}
 There is a coarsely surjective coarse map $\omega:\Z_+\to \Z^2$. Now $\soe(\omega):\soe(\Z_+)\to \soe(\Z^2)$ is not a surjective map obviously.
\end{ex}

\begin{lem}
\label{lem:soe1surjective}
 If $f:X\to Y$ is a coarse map between metric spaces, $Y$ is geodesic proper and $\soe(f):\soe(X)\to \soe(Y)$ is surjective then $f$ is already coarsely surjective.
\end{lem}
\begin{proof}
 Assume the opposite. Then $(\im f)^c\s Y$ contains a countable subset $(s_i)_i$ that is coarsely disjoint to $\im f$. Then by \cite[Proposition~93]{Hartmann2017a} there is a coarse ray $\rho:\Z_+\to Y$, an unbounded subsequence $(s_{i_k})_k$ and an entourage $E\s Y^2$ such that
\[
 (s_{i_k})_k\s E[\rho(\Z_+)]
\]
Now $\rho$ represents a point $r\in \soe(Y)$ and $\soe(f)$ is surjective. Thus there is some $p\in \soe(X)$ such that $\soe(f)(p)=r$. Suppose $p$ is represented by $\varphi:\Z_+\to X$ then $f\circ \varphi$ represents $r$ and has image in $\im f$. Thus there is an entourage $F\s Y^2$ such that $\rho(\Z_+)\s F[\im f]$. Then 
\[
 (s_{i_k})_k\s E\circ F[\im f]
\]
a contradiction to the assumption.
\end{proof}

\begin{rem}
Note that the coarse map 
\f{
\Z_+&\to \Z_+\\
n&\mapsto \lfloor \sqrt n\rfloor
}
is not coarsely injective. Since every map $\soe(\Z_+)\to \soe(\Z_+)$ is an isomorphism we cannot conclude that the functor $\soe(\cdot)$ reflects isomorphisms.
\end{rem}

\begin{lem}
If two subsets $U,V$ coarsely cover a metric space $X$ then
 \[
  \soe(U\cap V)=\soe(U)\cap \soe(V).
 \]
 \end{lem}
\begin{proof}
 The inclusion $\soe(U\cap V)\s \soe(U)\cap \soe(V)$ is obvious.
  
  We show the reverse inclusion: if $p\in \soe(U)\cap \soe(V)$ then it is represented by $\phi:\Z_+\to U$ in $\soe(U)$ and $\psi:\Z_+\to V$ in $\soe(V)$. Then there is an entourage $E\s X^2$ such that
\[
 E[\psi(\Z_+)]=\phi(\Z_+).
\]
 Denote by $F$ the set of indices $i\in\Z_+$ such that for each $j\in\Z_+$ the inclusion
\[
 (\phi(i),\psi(j))\s E\cap (V^c\times U^c).
\]
holds. Since $E\cap (V^c\times U^c)$ is bounded the set $F$ is finite. Now we construct a coarse map $\varphi:\Z_+\to U\cap V$: for every $i\in \Z_+\ohne F$ do:
\begin{enumerate}
\item if $\phi(i)\in V$ then define $\varphi(i):=\phi(i)$;
\item if $\phi(i)\not\in V$ then there exists some $j\in\Z_+$ such that $(\phi(i),\psi(j))\in E$. Now define $\varphi(i):=\psi(j)$.
\end{enumerate}
Fix a point $x_0\in U\cap V$ then for every $i\in F$ define: $\varphi(i)=x_0$. Then $\varphi$ represents $p$ in $\soe(U\cap V)$. 
\end{proof}

\begin{lem}
\label{lem:preservescoprod}
 The functor $\soe(\cdot)$ preserves finite coproducts.
\end{lem}
\begin{proof}
 Let $X=A\sqcup B$ be a coarse disjoint union of metric spaces. Without loss of generality we assume that $A,B$ cover $X$ as sets. Fix a point $x_0\in X$. Then there is a coarse map 
 \begin{align*}
 r:X&\to \Z\\
 x&\mapsto\begin{cases}
           d(x,x_0) & x\in A\\
           -d(x,x_0) & x\in B.
          \end{cases}
 \end{align*}
 Note that $\soe(\Z)=\{-1,1\}$ is a space which consists of two points with the discrete uniformity. Then $\soe(r)(\soe(A))=1$ and $\soe(r)(\soe(B))=-1$. Thus $\soe(X)$ is the uniform disjoint union of $\soe(A),\soe(B)$.
\end{proof}

\begin{prop}
\label{prop:hausdorff}
 Let $X$ be a metric space. The uniformity $\soe(X)$ is separated.
\end{prop}
\begin{proof}
If $p\not=q$ are two points in $\soe(X)$ we show there is a coarse cover $\ucover$ such that
\[
 q\not\in\ucover[p].
\]

Suppose $p$ is represented by $\phi:\Z_+\to X$ and $q$ is represented by $\psi:\Z_+\to X$. Now there is one of two cases:
 \begin{enumerate}
 \item there is a subsequence $(i_k)_k\s \N$ such that $\phi(i_k)_k\notclose \psi(\Z_+)$.  
 \item there is a subsequence $(j_k)_k\s \N$ such that $\psi(j_k)_k\notclose \phi(\Z_+)$.
 \end{enumerate}
 Without loss of generality we can assume the first case holds. By Lemma~\ref{lem:separationcover} there is a coarse cover $\ucover=\{U_1,U_2\}$ of $X$ such that $U_1\notclose \phi(i_k)_k$ and $U_2\notclose \psi(\Z_+)$. Then $q\not \in \ucover[p]$.

 Now 
 \f{
 q &\not\in\ucover[p]\\
 &=st(p,(\ucover[r])_r).
 }
 Thus the result.
\end{proof}

\begin{lem}
 \label{lem:totallybounded}
 If $X$ is a metric space, 
 \begin{itemize}
 \item $\ucover=(U_i)_{i\in I}$ is a coarse cover of $X$ and $p\in \soe(X)$ is represented by $\phi:\Z_+\to X$ then define
 \[
  I(p):=\{i\in I:\phi(\Z_+)\close U_i\}.
 \]
 If $S\s I$ is a subset then define
 \[
  U(S)=\{p\in \soe(X):\phi(\Z_+)\s E[\bigcup_{i\in S}U_i]\}
 \]
here $\phi:\Z_+\to X$ represents $p$ and $E\s X^2$ is an entourage. If $q\in \soe(X)$ then $q\in \ucover[p]$ if and only if $q\in U(I(p))$ and $p\in U(I(q))$.
\item Define
\[
 \ucover(S)=\{p\in \soe(X):p\in U(S),I(p)\z S\}.
\]
Then $q\in \ucover[p]$ if and only if there is some $S\s I$ such that $p,q\in \ucover(S)$. The uniform cover 
\[
(\ucover(S))_{S\s I}
\]
associated to $D_\ucover$ is a finite cover.
 \item The uniform space $\soe(X)$ is totally bounded.
 \end{itemize}
\end{lem}
\begin{proof}
\begin{itemize}
 \item easy.
 \item We just need to show: if $q\in \ucover[p]$ then $p\in \ucover(I(p)\cap I(q))$. For that it is sufficient to show if $\phi:\Z_+\to X$ represents $p$ then there is an entourage $E\s X^2$ such that $\phi(\Z_+)\s E[\bigcup_{i\in I(q)}U_i]$. Assume the opposite: there is some subsequence $(i_k)_k\s \Z_+$ such that
 \[
  \phi(i_k)_k\notclose \bigcup_{i\in I(p)\cap I(q)}U_i.
 \]
 Now $\phi(\Z_+)\notclose \bigcup_{i\not\in I(p)}U_i$ thus
 \[
  \phi(i_k)_k\notclose (\bigcup_{i\in I(p)\cap I(q)}U_i)\cup(\bigcup_{i\not\in I(p)}U_i).
 \]
And thus $ \phi(i_k)_k\notclose \bigcup_{i\in I(q)}U_i$ a contradiction to the assumption.
 \item easy
\end{itemize}
\end{proof}

\begin{notat}
 If $A,B\s X$ are two subsets of a metric space and $x_0\in X$ a point then define
 \f{
  \chi_{A,B}:\N&\to \R_+\\
  i&\mapsto d(A\ohne B(x_0,i),B\ohne B(x_0,i))
 }
Now $A\close B$ if and only if $\chi_{A,B}$ is bounded. And if $A_1,A_2,B\s X$ are subsets with $A_1=E[A_2]$ then the function
 \[
  i\mapsto|\chi_{A_1,B}(i)-\chi_{A_2,B}(i)|
 \]
 is bounded.
 \end{notat}

\begin{defn}
\label{defn:metricone}
 Let $X$ be a metric space. If two endpoints $p,q\in \soe(X)$ are represented by coarse maps $\phi,\psi:\Z_+\to X$ then the distance of $p$ to $q$ is at least $f\in\R_+^\N$, written $d(p,q)\ge f$, if there is a subsequence $(i_k)_k\s \Z_+$ such that one of the following holds
 \begin{enumerate}
  \item $\phi(\Z_+)\notclose \psi(i_k)_k$ and $\chi_{\phi(\Z_+),\psi(i_k)_k}+c\ge f$ for some $c\ge 0$;
  \item $\psi(\Z_+)\notclose \phi(i_k)_k$ and $\chi_{\psi(\Z_+),\phi(i_k)_k}+c\ge f$ for some $c\ge 0$.
 \end{enumerate}
 \item We define $d(p,q)=0$ if and only if $p=q$.
\end{defn}

\begin{lem}
\label{lem:chiprops}
If $X$ is a metric space and $\ucover$ a coarse cover of $X$ then there is an unbounded function $f\in \R_+^\N$ such that for every two endpoints $p,q\in \soe(X)$ the relation $q\not\in \ucover[p]$ implies $d(p,q)\ge f$.
\end{lem}
\begin{proof}
 By Lemma~\ref{lem:totallybounded} the uniform space $\soe(X)$ is totally bounded. Without loss of generality we can fix an endpoint $p\in \soe(X)$ and study the endpoints $q\in \soe(X)$ for which $q\not\in \ucover[p]$.

 We will define a function $f\in \R_+^\N$ as the minimum of a finite collection $f_0,\ldots,f_n\ge 0$ of numbers. 

 If $q\not\in \ucover[p]$ there are 2 cases:
 \begin{enumerate}
 \item $p\not\in U(I(q))$: There is a subset $(i_k)_k\s \Z_+$ such that 
 \[
 \phi(i_k)_k\notclose \bigcup_{i\in I(q)}U_i.
 \]
 Now for every $S\s I$ if $p\not\in U(S)$: then there is some subset $(i_k)_k\s \Z_+$ such that $\phi(i_k)_k\notclose \bigcup_{i\in S}U_i$. Thus if $I(q)= S$ then $d(p,q)> \chi_{\phi(i_k)_k,\bigcup_{i\in S}U_i}$. Define
 \[
  f_S:=\chi_{\phi(i_k)_k,\bigcup_{i\in S}U_i}
 \]
in this case.
 \item $q\not\in U(I(p))$: There is some $(i_k)_k\s \Z_+$ such that 
 \[
  \psi(i_k)_k\notclose \bigcup_{i\in I(p)}U_i.
 \]
 Now $\psi(i_k)_k\s \bigcup_{i\not\in I(p)}U_i$ and $\phi(\Z_+)\notclose \bigcup_{i\not\in I(p)}U_i$. Then define
 \[
  f_b:=\chi_{\phi(\Z_+),\bigcup_{i\not\in I(p)}U_i}.
 \]
 \end{enumerate}
 Now define
 \[
  f(i):=\min_S(f_S(i),f_b(i))
 \]
 for every $i$. Then $f$ has the desired properties.
\end{proof}

\begin{prop}
\label{prop:countablebase}
 If $X$ is a metric space then the uniformity on $\soe(X)$ is coarser than the uniformity $\tau_d$, induced by $d$.
\end{prop}
\begin{proof}
 By Lemma~\ref{lem:chiprops} every entourage in $\soe(X)$ is a neighborhood of an entourage of $\tau_d$ on $\soe(X)$. 
 \end{proof}

\section{Side Notes}
\label{sec:side}

\begin{lem}\name{Higson corona}
 If $X$ is a metric space then the $C^*$-algebra that determines the Higson corona is a sheaf. That means exactly that 
the association
 \[
  U\mapsto C(\nu U)=B_h(U)/B_0(U)
 \]
for every subset $U\s X$ is a sheaf with values in $\cstar$. By a sheaf we mean a sheaf on the Grothendieck topology determined by coarse covers on subsets of a coarse space.
\end{lem}
\begin{proof}
 We recall a few definitions which can be found in \cite[p.29, 30]{Roe2003}.
 \begin{itemize}
  \item The algebra of bounded functions that satisfy the Higson condition is denoted by $B_h$.
  \item A bounded function $f:X\to \C$ satisfies the Higson condition if for every entourage $E\s X^2$ the function
  \f{
  df|_E:E&\to \C\\
  (x,y)&\mapsto f(y)-f(x)
  }
  tends to $0$ at infinity.
  \item the ideal of bounded functions that tend to $0$ at infinity is called $B_0$.
  \item A function $f:X\to \C$ tends to $0$ at infinity if for every $\varepsilon>0$ there is a bounded subset $B\s X$ 
such that $|f(x)|\ge \epsilon$ implies $x\in B$.
 \end{itemize}
We check the sheaf axioms:
\begin{enumerate}
 \item global axiom: if $U_1,U_2$ coarsely cover a subset $U\s X$ and $f\in B_h(X)$ such that $f|_{U_1}\in B_0(U_1)$ 
and $f|_{U_2}\in B_0(U_2)$ we show that $f\in B_0(U)$ already. Let $\varepsilon>0$ be a number. Then there are bounded 
subsets $B_1\s U_1$ and $B_2\s U_2$ such that $|f|_{U_i}(x)|\ge \varepsilon$ implies $x\in B_i$ for $i=1,2$. 
Now
\[
 B:=B_1\cup B_2 \cup (U_1\cup U_2)^c
\]
is a bounded subset of $U$. Then $|f(x)|\ge \varepsilon$ implies $x\in B$. Thus $f\in B_0(U)$.
\item gluing axiom: if $U_1,U_2$ coarsely cover a subset $U\s X$ and $f_1\in B_h(U_1),f_2\in B_h(U_2)$ are functions such that
\[
 f_1|_{U_2}=f_2|_{U_1}+g
\]
where $g\in B_0(U_1\cap U_2)$. We show there is a function $f\in B_h(U)$ which restricts to $f_1$ on $U_1$ and $f_2+g$ 
on $U_2$. Define:
\begin{align*}
 f:U&\to \C\\
 x&\mapsto\begin{cases}
           f_1(x) & x\in U_1\\
           f_2(x)+g & x\in U_2\\
           0 & \mbox{otherwise}
          \end{cases}
\end{align*}
then $f$ is a bounded function. We show $f$ satisfies the Higson condition: Let $E\s U^2$ be an entourage and
$\varepsilon>0$ be a number. Then there are bounded subsets $B_1\s U_1$ and $B_2\s U_2$ such that $|df_i|_{E\cap 
U_i^2}(x,y)|\ge \varepsilon$ implies $x\in B_i$ for $i=1,2$. There is a bounded subset $B_3\s U$ such that
\[
 E\cap (U_1^2\cup U_2^2)^c\cap U^2\s B_3^2.
\]
Define
\[
 B:=B_1\cup B_2\cup B_3
\]
then $|df|_E(x)|\ge\varepsilon$ implies $x\in B$. Thus $f$ has the desired properties.
\end{enumerate}
\end{proof}

\begin{prop}
If $X$ is a proper geodesic metric space denote by $\sim$ the relation on $\soe(X)$ of belonging to the same uniform connection component in $\soe(X)$. Then there is a continuous bijection
 \[
  \soe(X)/\sim\to \Omega(X)
 \]
where the right side denotes the space of ends of $X$ as a topological space.
\end{prop}
\begin{proof}
There are several different definitions for the space of ends of a topological space. We use~\cite[Definition~8.27]{Bridson1999}.

 An end in $X$ is represented by a proper continuous map $r:[0,\infty)\to X$. Two such maps $r_1,r_2$ represent the same end if for every compact subset $C\s X$ there is some $N\in\N$ such that $r_1[N,\infty),r_2[N,\infty)$ are contained in the same path component of $X\ohne C$.

 If $r:[0,\infty)\to X$ is an end then there is a coarse map $\varphi:\Z_+\to X$ and an entourage $E\s X^2$ such that
 \[
  E[r[0,\infty)]=\varphi(\Z_+).
 \]
We construct $\varphi$ inductively:
\begin{enumerate}
 \item $\varphi(0):=r(0)$
 \item if $\varphi(i-1)=r(t_{i-1})$ is already defined then $t_i:=\min\{t>t_{i-1}:d(r(t_{i-1}),r(t))=1\}$. Set $\varphi(i):=r(t_i)$.
\end{enumerate}
By the above construction $\varphi$ is coarsely uniform. The map $\varphi$ is coarsely proper because $r$ is proper and $X$ is proper.

Note that every geodesic space is also a length space. If for some compact subset $C\s X$ the space $X\ohne C$ has two path components $X_1,X_2$ then for every $x_1\in X_1,x_2\in X_2$ a path (in particular the shortest) joining $x_1$ to $x_2$ contains a point $c\in C$. Thus
\[
 d(x_1,x_2)=\inf_{c\in C}(d(x_1,c)+d(x_2,c))
\]
Then $X$ is the coarse disjoint union of $X_1,X_2$. On the other hand if $X$ is the coarse disjoint union of subspaces $X_1,X_2$ then there is a bounded and in particular because $X$ is proper compact subset $C\s X$ such that
\[
 X\ohne C=X_1'\sqcup X_2'
\]
is a path disjoint union and $X_1'\s X_1,X_2'\s X_2$ differ only by bounded sets.

Now we show the association is continuous:

We use \cite[Lemma~8.28]{Bridson1999} in which $\mathcal G_{x_0}(X)$ denotes the set of geodesic rays issuing from $x_0\in X$. Then \cite[Lemma~8.28]{Bridson1999} states that the canonical map 
 \[
  \mathcal G_{x_0}\to \Omega(X)
 \]
 is surjective. Fix $r\in \mathcal G_{x_0}$. Then $\tilde V_n\s \mathcal G_{x_0}$ denotes the set of proper rays $r':\R_+\to X$ such that $r'(n,\infty),r(n,\infty)$ lie in the same path component of $X\ohne B(x_0,n)$. Now \cite[Lemma~8.28]{Bridson1999} states the sets $(V_n=\{[r']:r'\in \tilde V_n\})_n$ form a neighborhood base for $[r]\in \Omega(X)$.

 Now to every $n$ we denote by $U_1^n$ the path component of $X\ohne B(x_0,n)$ that contains $r(\R_+)$ and we define $U_2^n:=X\ohne U_1^n$. For every $n\in\N$ the sets $U_1^n,U_2^n$ are a coarse cover of $X$.
 
 Suppose $\rho:\Z_+\to X$ is a coarse map associated to $r$ and represents $\tilde r\in \soe(X)$. If $s\in \mathcal G_0$  suppose $\sigma:\Z_+\to X$ is the coarse map associated to $s$ and represents $\tilde s\in \soe(X)$. If $[s]\not\in V_n$ then $\sigma(\Z_+)\notclose U_1^n$. This implies $\tilde s\not\in \{U_1^n,U_2^n\}[\tilde r]$.
 
 Thus for every $n\in \N$ there is an inclusion $\{U_1^n,U_2^n\}[\tilde r]/\sim\s V_n$ by the association.
\end{proof}

\section{Remarks}

The starting point of this research was an observation in the studies of~\cite{Hartmann2017a}: coarse cohomology with twisted coefficients looked like singular cohomology on some kind of boundary. We tried to find a functor from the coarse category to the category of topological spaces that would reflect that observation.

And then we noticed that two concepts play an important role: One is the choice of topology on the \emph{space of ends} and one is the choice of points. The points were designed such that 
\begin{itemize}
 \item coarse maps are mapped by the functor to maps of sets 
 \item and the space $\Z_+$ is mapped to a point
\end{itemize}
If the metric space is Gromov hyperbolic then coarse rays represent the points of the Gromov boundary, thus the Gromov boundary is a subset of the space of ends. The topology was trickier to find. We looked for the following properties:
 \begin{itemize}
  \item coarse maps are mapped to continuous maps
  \item coarse embeddings are mapped to topological embeddings
 \end{itemize}
 Now a proximity relation on subsets of a topological space helps constructing the topology on the space of ends of Freudenthal. We discovered that coarse covers on metric spaces give rise to a totally bounded uniformity and thus used that a uniformity on a space gives rise to a topology.

Finally, after a lucky guess, we came up with the uniformity on the set of endpoints. In which way does the space of ends functor reflect isomorphism classes will be studied in a paper that follows.

It would be possible, conversely, after a more thorough examination to find more applications. Coarse properties on metric spaces may give rise to topological properties on metrizable uniform spaces.

We wonder if this result will be of any help with classifying coarse spaces up to coarse equivalence. However, as of yet, the duality has not been studied in that much detail.

\bibliographystyle{utcaps}
\bibliography{mybib}

\address

\end{document}